\documentclass[onefignum,onetabnum]{siamart190516}

\usepackage{lipsum}
\usepackage{amsfonts}
\usepackage{graphicx}
\usepackage{epstopdf}
\usepackage[caption=false]{subfig} %
\usepackage{algorithmic}

\usepackage{hyperref}
\usepackage{cleveref}

\usepackage{tikz}
\usetikzlibrary{arrows}
\usepackage{pgfplots}
\pgfplotsset{width=7cm, compat=1.10}
\usepgfplotslibrary{fillbetween}
\tikzstyle{block}=[draw opacity=0.7,line width=1.4cm]
\usetikzlibrary{intersections, calc,angles,quotes,babel,matrix, arrows,shadows,patterns,automata}

\ifpdf
\DeclareGraphicsExtensions{.eps,.pdf,.png,.jpg}
\else
\DeclareGraphicsExtensions{.png}
\fi

\usepackage{chngcntr}
\counterwithout{equation}{section}
\renewcommand{\theequation}{\arabic{equation}}

\newsiamremark{remark}{Remark}
\newsiamremark{hypothesis}{Hypothesis}
\newsiamthm{claim}{Claim}

\newcommand{\NN}{\mathbb{N}}
\newcommand{\R}{\mathbb{R}}

\headers{Aggregation-Diffusion Equations for Collective Behaviour}{R. Bailo, J.A. Carrillo, D. G\'omez-Castro}

\title{Aggregation-Diffusion Equations for Collective Behaviour in the Sciences\thanks{The authors were supported by the Advanced Grant Nonlocal-CPD (Nonlocal PDEs for Complex Particle Dynamics: Phase Transitions, Patterns and Synchronization) of the European Research Council Executive Agency (ERC) under the European Union's Horizon 2020 research and innovation programme (grant agreement No. 883363). JAC is partially supported by the “Maria de Maeztu” Excellence
Unit IMAG, reference CEX2020-001105-M, funded by
MCIN/AEI/10.13039/501100011033/.
DGC is support by RYC2022-037317-I from the Spanish government.
}}

\author{
Rafael Bailo%
\thanks{Mathematical Institute, University of Oxford, Oxford, OX2 6GG, UK. 
(\email{bailo@maths.ox.ac.uk})} 
\and 
Jos\'e A. Carrillo%
\thanks{Mathematical Institute, University of Oxford, Oxford, OX2 6GG, UK. (\email{carrillo@maths.ox.ac.uk}%
)}
\and 
David G\'omez-Castro%
\thanks{Departamento de Matemáticas, Universidad Autónoma de Madrid, 28049 Madrid, Spain (\email{david.gomezcastro@uam.es}})
}

\usepackage{amsopn}

\DeclareMathOperator*{\argmin}{argmin}

\ifpdf
\hypersetup{
pdftitle={An Example Article},
pdfauthor={D. Doe, P. T. Frank, and J. E. Smith}
}
\fi

\newcommand{\RB}[1]{{\color{blue} ([RB] #1)}}
\newcommand\rb\RB

\usepackage{mathtools}
\usepackage{suffix}

\newcommand{\nhalf}{\frac{1}{2}}

\renewcommand{\i}{_{i}}
\newcommand{\ip}{_{i+1}}

\newcommand{\ih}{_{i+\nhalf}}
\newcommand{\imh}{_{i-\nhalf}}

\newcommand{\n}{^{n}}
\newcommand{\np}{^{n+1}}

\newcommand{\nss}{^{**}}

\newcommand{\Dt}{\Delta t}
\newcommand{\Dx}{\Delta x}

\DeclarePairedDelimiter{\prt}{(}{)}

\DeclarePairedDelimiter{\abs}{|}{|}
\DeclarePairedDelimiter{\norm}{\|}{\|}
\DeclarePairedDelimiter{\set}{\{}{\}}

\DeclarePairedDelimiter{\positive}{(}{)^{+}}
\DeclarePairedDelimiter{\negative}{(}{)^{-}}

\newcommand\pos\positive
\renewcommand\neg\negative

\WithSuffix\newcommand\pos*{\positive*}
\WithSuffix\newcommand\neg*{\negative*}

\newcommand{\Wik}{W_{i-k}}

\renewcommand{\k}{_{k}}

\renewcommand{\d}{\mathrm{d}}
\newcommand{\dd}{\mathop{}\!\d}

\newcommand{\der}[2]{\frac{\d #1}{\d #2}}
\newcommand{\dx}{\dd x}
\newcommand{\dy}{\dd y}

\begin{document}

\maketitle

\begin{abstract}
This is a survey article based on the content of the plenary lecture given by Jos\'e A. Carrillo at the ICIAM23 conference in Tokyo. It is devoted to produce a snapshot of the state of the art in the analysis, numerical analysis, simulation, and applications of the vast area of aggregation-diffusion equations. We also discuss the implications in mathematical biology explaining cell sorting in tissue growth as an example of this modelling framework. This modelling strategy is quite successful in other timely applications such as global optimisation, parameter estimation and machine learning.
\end{abstract}

\section{Introduction}

A surprisingly large number of complex phenomena can be modelled as systems of point particles that interact under short and long-range forces. This modelling framework, which is classical in the study of many-body systems (e.g. gravitational collapse, electron transport in semi-conductors, or kinetic equations in statistical physics), has recently found applications in other areas of physics, as well as in biology, social sciences, machine learning, and even optimisation, see for instance \cite{Carillo09_Kinetic_Attraction-Repulsion,carrillo2019population,CHSV,CCTT,BBG2024}, and is a prolific and current field of research. The applications are diverse, and span both the microscopic and macroscopic scales: ion channel transport, chemotaxis, bacterial orientation, cellular adhesion, angiogenesis, animal herding, or human crowd motion are but a few examples. The interaction effects between the ``individuals'' in each example can be approximated by long-range attractive forces (e.g. ligand binding, electrical interactions, or social preferences) and short-range repulsion effects (due to volume constraints or crowding).

\begin{figure}
	\centering
	\includegraphics[width=0.8\textwidth,clip]{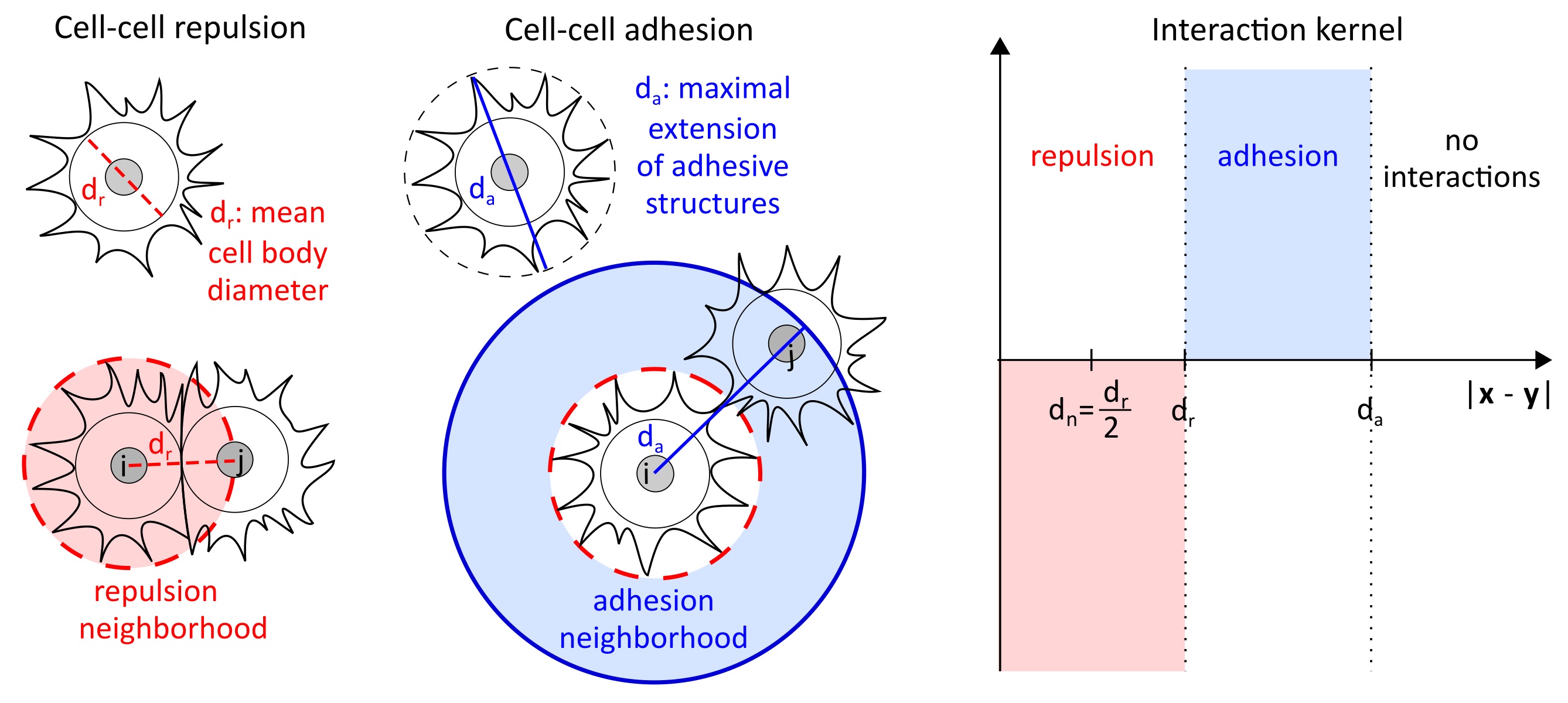}\\
	\caption{Interaction forces of an agent-based model for cell adhesion, see \cite{CARRILLO201875,carrillo2019population}.}
	\label{fig:Cellinteractions}
\end{figure}

There are several agent-based (microscopic) modelling approaches to describe these phenomena, such as cellular automata or driven Brownian particles, all capable of encoding complex behaviour. Yet, said behaviour is best understood and rigorously analysed by studying the mean-field limit (or its variants): the limit of the microscopic model as the number of agents becomes large, see \cite{oelschlager1990large,BCC11,J}.
As an example, suppose that we study two types of cells, each able to express different surface proteins and ligands (such us cadherins or nectins). Assume, for cells of the first kind, that nuclei are located at positions $\{ x_i \}$, $i=1,\dots,N$, and, for the second, at $\{ y_j \}$, $j=1,\dots,M$; let also $N=M$ for simplicity. We consider cells that interact either by attraction at medium distances (e.g. through filopodia, protrusions of the cellular membrane), or by strong repulsion at short distances (due to volume and size constraints); a diagram is given in \cref{fig:Cellinteractions}. If the force exerted by cells of type $b$ onto type $a$ is radial and conservative (with potential $W^N_{ab}$, where $a,b=1,2$), we may pose the elementary agent-based model
\begin{align*}%
	\dot x_i & = - \frac{m_1}N \sum_{j\neq i}\nabla W^N_{11}(x_i-x_j) - \frac{m_2}N\sum_{j} \nabla W^N_{12}(x_i-y_j),              \\
	\dot y_i & = - \frac{m_2}N \sum_{j\neq i} \nabla W^N_{22}(y_i-y_j) - \frac{m_1}N\sum_{j} \nabla W^N_{21}(y_i-x_j),\label{aux2}
\end{align*}
where $m_1$ and $m_2$ are the typical cellular masses of each cell kind.

Under suitable assumptions, the empirical measures (weighted sums of Dirac deltas supported at each particle's position) associated to the model approximate the macroscopic normalised cell densities in the many-particle limit:
\begin{align*}
	\rho_1(t,x) \simeq \frac{m_1}N \sum_{i=1}^N \delta_{x_i(t)}{(x)} \quad \textrm{and} \quad \rho_2(t,x) \simeq \frac{m_2}N \sum_{i=1}^N \delta_{y_i(t)}{(x)} \quad \textrm{as } N\to \infty.
\end{align*}
When considering the biological context, it is reasonable to suppose that the attraction forces act only within a cut-off radius $R$, and that the repulsion is purely localised \cite{CalvezCarrillo}; a natural scaling choice for the potentials is therefore $W^N_{ab}\simeq \varepsilon \delta_0 + W_{ab}$ as $N\to \infty$ \cite{Oes,carrillo2019population}, where $\varepsilon$ relates to the typical volume occupied by each cell (assumed here equal across species and interactions). In this scaling, the limit $N\to\infty$ leads to the mean-field model
\begin{equation}
	\label{model_chs}
	\left\{
	\begin{aligned}
		\partial_t \rho_1 & = \nabla \cdot\Big(\rho_1 \nabla \big(\varepsilon (\rho_1+\rho_2) + W_{11}\ast\rho_1 + W_{12}\ast\rho_2 \big)\Big), \\
		\partial_t \rho_2 & = \nabla \cdot\Big(\rho_2 \nabla \big(\varepsilon (\rho_1+\rho_2) + W_{22}\ast\rho_2 + W_{21}\ast\rho_1 \big)\Big).
	\end{aligned}
	\right.
\end{equation}

Each of the equations in the system above is a particular case of the \emph{aggregation-diffusion equation}
\begin{align} \label{aggeqn0}
	\frac{\partial\rho}{\partial t}+\nabla \cdot \left(\rho u\right)=0 ,
	\quad u=-\nabla \xi,
	\quad \xi=U' \left ( \rho \right )+ V + W*\rho ,
\end{align}
a generic continuum equation for the kinematic evolution of a population density $\rho(t,x)$. The velocity field $u$ drives the evolution of the density through a balance of forces: diffusion, resulting from localised repulsion effects or noise, which is modelled by the convex function $U(\rho)$; drift, due to external forces, described by the potential $V(x)$, and a symmetric interaction force (aggregation, repulsion, or both), given by the potential $W(x)$, see \cite{BodnarVelazquez,TopazBertozzi2}. The diffusion is commonly given by $U(s)=s\log s$ (linear diffusion) or $U(s)={s^m}/{(m-1)}$ for $m>0$ (non-linear diffusion \cite{Va}). In these cases, the diffusion coefficient is $m \rho^{m-1}$: in the \emph{porous medium} range $m > 1$, diffusion is stronger where the concentration is large; in the \emph{fast diffusion} regime, $m < 1$, where the concentration is small. So-called confining potentials, $V = |x|^p$ for $p > 0$, often appear in Fokker-Planck or McKean-Vlasov equations \cite{cmcv-03,carrillo2020LongTimeBehaviourPhase,CarrilloGomezCastroVazquez2022JMPA}.

\section{Gradient Flows}\label{sec:gradient_flows}

\begin{figure}
	\centering
	\includegraphics[width=0.2\textwidth,trim={3cm 2cm 2.5cm 1.25cm},clip]{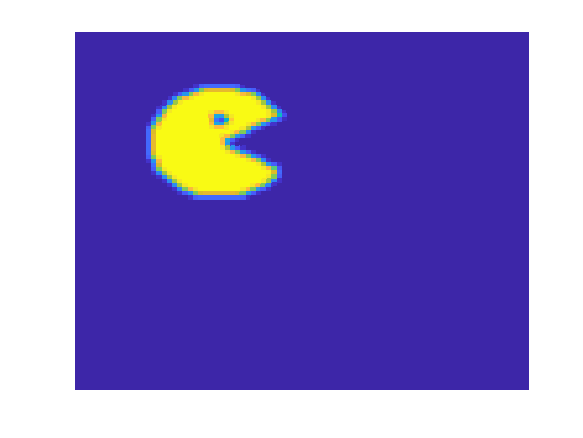}
	\includegraphics[width=0.2\textwidth,trim={3cm 2cm 2.5cm 1.25cm},clip]{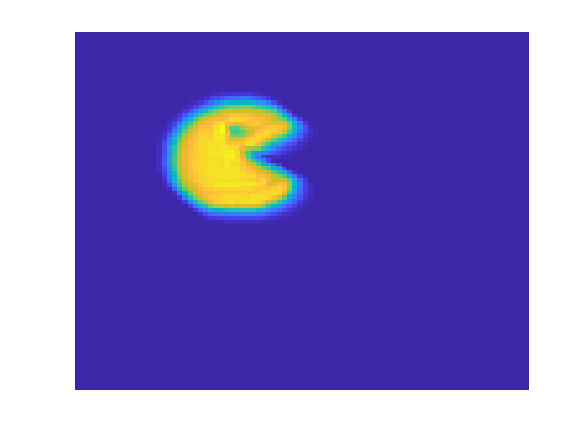}
	\includegraphics[width=0.2\textwidth,trim={3cm 2cm 2.5cm 1.25cm},clip]{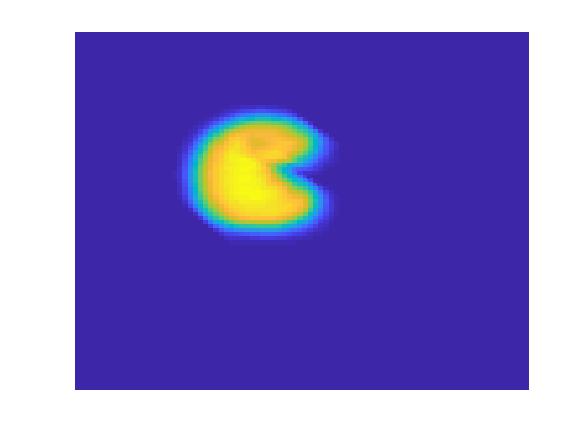}\\
	\includegraphics[width=0.2\textwidth,trim={3cm 2cm 2.5cm 1.25cm},clip]{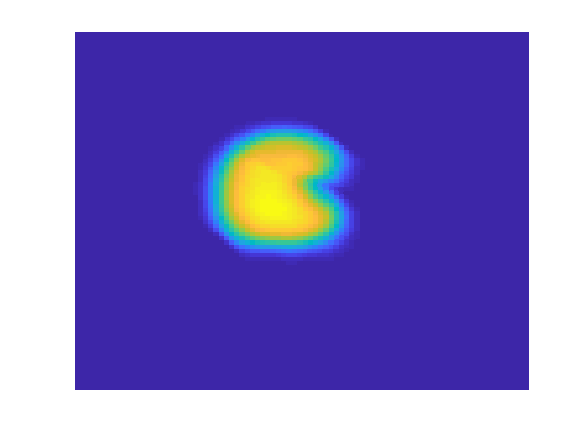}
	\includegraphics[width=0.2\textwidth,trim={3cm 2cm 2.5cm 1.25cm},clip]{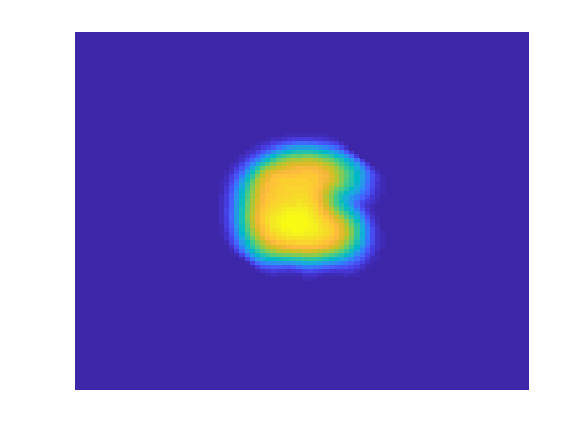}
	\includegraphics[width=0.2\textwidth,trim={3cm 2cm 2.5cm 1.25cm},clip]{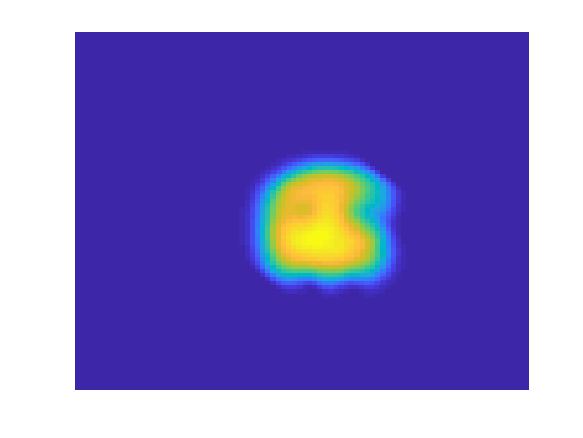}\\
	\includegraphics[width=0.2\textwidth,trim={3cm 2cm 2.5cm 1.25cm},clip]{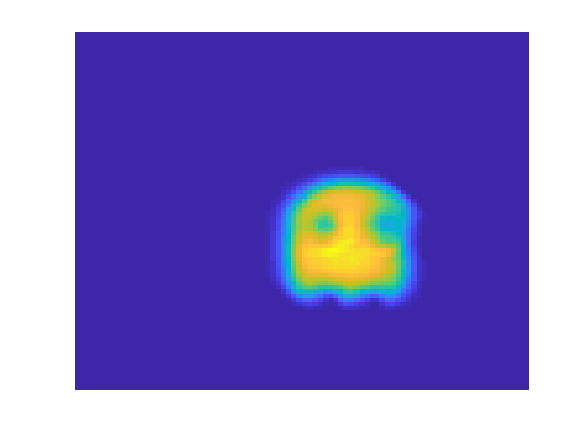}
	\includegraphics[width=0.2\textwidth,trim={3cm 2cm 2.5cm 1.25cm},clip]{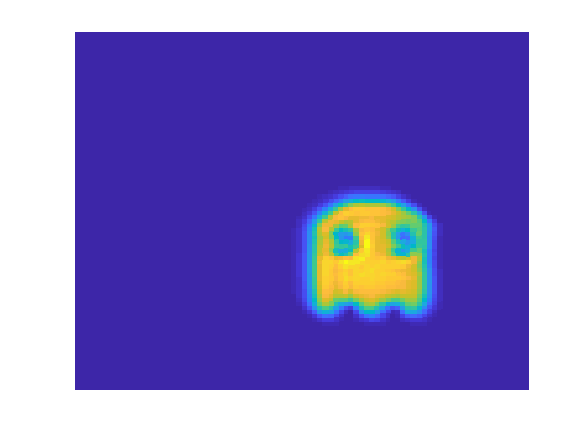}
	\includegraphics[width=0.2\textwidth,trim={3cm 2cm 2.5cm 1.25cm},clip]{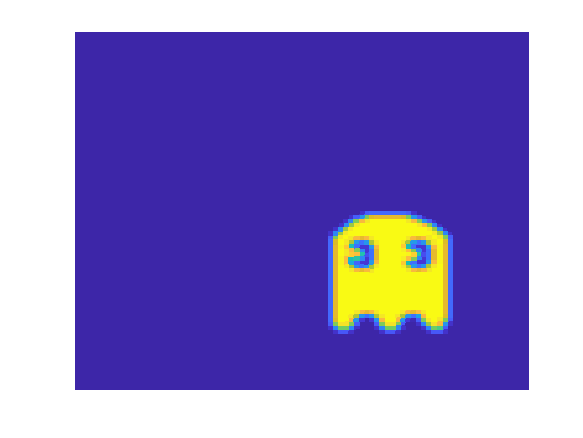}\\
	\caption{Numerical approximation of a geodesic curve using the $2$-Wasserstein distance between the characteristic sets of Pac-Man and the Ghost (suitably normalised). \href{https://figshare.com/articles/media/Wasserstein_Geodesic_between_PacMan_and_Ghost/7665377}{Video available online}.}\label{fig:pacman}
\end{figure}

The aggregation-diffusion equation \eqref{aggeqn0} can be understood as the gradient flow of the free-energy functional
\begin{equation}\label{free}
	\mathcal F[ \rho ]
	= \int_{\R^d} U(\rho (x)) \dx
	+ \int_{\R^d} V(x) \rho(x) \dx
	+ \frac{1}{2}\int_{\R^d}\int_{\R^d} \rho(x) W(x-y) \rho(y) \dy\dx ,
\end{equation}
in a suitable metric space: the $2$-Wasserstein space of probability densities over $\mathbb R^d$, see \cite{JKO,Villani,AGS}. This space arises by equipping the set of probability densities with finite second moment, denoted $\mathcal P_2 (\mathbb R^d)$, with the $2$-Wassertein distance, denoted $d_2$. Said distance, which originated in the theory of optimal transport, has found fruitful applications in the study of mean-field limits because it induces a topology where the evolution of empirical measures in time is continuous, unlike in Sobolev spaces. As a result, $d_2$ is helpful when quantifying the convergence of particle systems to their mean-field limit. This notion of distance for probabilities may also be used to define transport-based interpolations between densities as geodesic curves in the $2$-Wasserstein space; \cref{fig:pacman} shows the geodesic interpolation between the normalised characteristics of two sets.

Due to the lack of linear structure in the $2$-Wassertein space, the precise construction of solutions to the aggregation-diffusion equation \eqref{aggeqn0} as gradient flows of the free energy functional \eqref{free} is quite involved. Roughly, a discrete-in-time sequence of measures is generated from the datum $\rho_0\in{\mathcal P}_2(\R^d)$ via the so-called \emph{JKO iteration}
\begin{align*}
	\rho_{k+1}= \argmin_{\rho\in {\mathcal P}_2(\R^d)} \left\{ \frac{1}{2 \Delta t}d_2^2(\rho,\rho_{k})+\mathcal F[\rho] \right\},
\end{align*}
for $k\in\NN$, and for any fixed $\Delta t>0$ \cite{JKO}. A suitable time interpolation leads to a curve of measures $\rho_{\Delta t}$ that converge to the unique solution of \eqref{aggeqn0} as $\Delta t\to 0$, see \cite{JKO,AGS}. Solutions constructed by this variational scheme satisfy the fundamental \emph{energy dissipation identity}:
\begin{equation}\label{dissipation}
	\der{}{t} \mathcal F[ \rho ]
	= \int_{\R^d} \xi \partial_t \rho \dx
	= \int_{\R^d} \xi \nabla \cdot \left(\rho \nabla \xi \right) \dx
	= - \int_{\R^d} \rho \abs*{\nabla \xi}^2 \dx
	\leq 0 \,,
\end{equation}
under certain assumptions on the potentials and the internal energy.

If the functional $\mathcal F[\rho]$ has a unique global minimiser in $\mathcal P_2 (\mathbb R^d)$, then the solution to \eqref{aggeqn0} converges in time to the minimiser whenever the functional satisfies a suitable notion of convexity (uniform convexity along the geodesics induced by $d_2$). However, this may not be the case; the purely diffusive setting ($V=W=0$) leads to asymptotic spreading of $\rho(t)$ in self-similar profiles for large times \cite{Va}. In general, the asymptotic balance between aggregation and diffusion is quite subtle, and significant effort has been devoted to its study. Under certain conditions, the existence of radially decreasing, compactly supported, global minimisers of $\mathcal F[\rho]$ (unique up to translations) has been established \cite{CCP15,Carrillo2019a}. In some extreme cases, the minimisation of $\mathcal F$ leads to the formation of Dirac deltas \cite{CDFL} that accumulate only a fraction of the total mass. When it comes to the dynamical problem, understanding these features is a considerable challenge, see the recent surveys \cite{CCY,DavidSEMA}. Almost-sharp conditions on $W$ have been given in \cite{Carrillo+GC+Yao+Zheng2021} such that $\rho(t)$ behaves like the solution of the purely diffusive problem. Yet, in some settings \cite{CarrilloGomezCastroVazquez2022JMPA}, the solution can develop asymptotically the Dirac delta suggested by the minimisation of the free energy.

\section{Applications in Mathematical Biology}\label{sec:math_bio}

\begin{figure}
	\centering
	\includegraphics[width=1.0\linewidth]{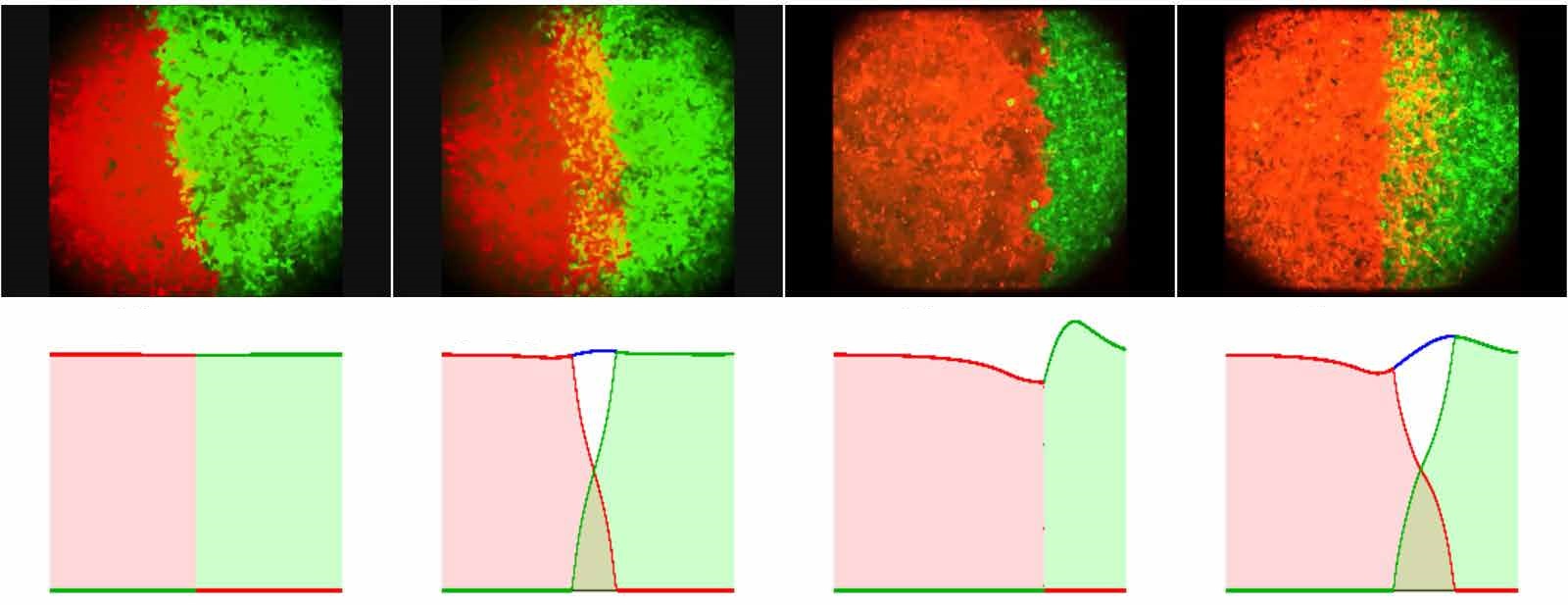}
	\caption{\emph{In vitro} experiments from \cite{katsunuma2016synergistic} compared to \emph{in silico} experiments from \cite{carrillo2019population}. The numerical schemes are based on finite volumes, as in \cite{BCH2020,BCH2023}. \href{https://figshare.com/articles/media/Front_propagation_and_intermingling_of_cell_types_experiments_versus_mathematical_model_simulations_/7707890}{Video available online}.}
	\label{fig:figtethumbnail}
\end{figure}

Mathematical biology has a special place among the applications of aggregation-diffusion systems. Cell population models are particularly relevant, and often use the aggregation-diffusion equation, either in scalar or in system form, to describe the motion of cells as they concentrate or separate from a target. The diffusion terms described above are consistent with population pressure effects, whereby groups of cells naturally spread away from areas of high concentration. The aggregation term has been classically used as a model for chemotaxis, as cells often direct their movement along the increasing gradient of a chemical agent produced either by external sources or by the very cells. The by-now classical Keller-Segel model for chemotaxis \cite{KellerSegel} is an archetypical example of the aggregation-diffusion equation in mathematical biology.

A new paradigm in the modelling of cell adhesion suggests that ligand binding through filopodia can be captured by short range interaction potentials, as suggested in \cref{fig:Cellinteractions}. An example of such model was developed in \cite{carrillo2019population}: it considers two species of cells, each attracted to itself and to the other in the form described by \eqref{model_chs}. Despite its apparent simplicity, the model can describe several realistic behaviours. If the all the attraction strengths are equal, all cells simply mix evenly. However, if the relative strengths vary, intricate phenomena arise; for instance, one species of cells might completely envelop the other, forming a ``halo''. \cref{fig:figtethumbnail} shows the astonishing match between the solutions of the model and the experimental results of \cite{katsunuma2016synergistic}. This model demonstrates that the differential adhesion hypothesis introduced by Steinberg \cite{steinberg1962I} is the driving mechanism for cell sorting \cite{doi:10.1137/22M1506079}, modelled here through the asymmetry on the interaction potentials between the two type of cells. Aggregation-diffusion systems \eqref{model_chs} with such asymmetry do not typically have a gradient-flow structure, and can develop a wealth of complex phenomena: among others, travelling wave solutions, free boundaries, and mixing or segregation between species; each and every one of them a challenging open problem in the analysis of Partial Differential Equations.

Beyond biology, the applications of these ideas in global optimisation and machine learning are not only timely, but potentially quite impactful. Many researchers are exploring them in order to perform gradient-free optimisation \cite{PTTM17,CHSV,CJLZ21,BBG2024,carrillo2024interacting}, gradient-free sampling \cite{GHLS20,korba2020non,DL21,CV21,CHSV}, Bayesian parameter estimation \cite{BFHM17,LZTM19,GSW19,LMT22,LL22,carrillo2024sparse}, or neural network calibration \cite{mokrov2021large,chizat2022infinite,fernandez2022continuous,lee2023deep}.

\section{Analytical Results}

There is a rich theory of the aggregation equation \eqref{aggeqn0}, yet its analysis remains a very active field, producing an array of results each year. This section will provide an overview, referring the reader to the recent surveys \cite{CCY,DavidSEMA} for detailed explanations and additional references.

\subsection{Local Minimisers of the Free Energy}
\label{sec:Euler Lagrange}
The Euler-Lagrange equation for the minimisation of the free energy \eqref{free} in ${\mathcal P}_2(\R^d)$ is
\begin{equation*}
	\begin{dcases}
		U'(\rho) + V + W*\rho = C   & \text{if } \rho > 0,    \\
		U'(\rho) + V + W*\rho \ge C & \text{in } \mathbb R^d.
	\end{dcases}
\end{equation*}
The constant $C$ can sometimes be expressed in terms of the energy of candidates for the minimiser of the free energy. Global minimisers may have a multiply connected support, yet $C$ must be the same on each, which is not the case for local minimisers. If the map $U' : (0,\infty) \to (0,\infty)$ is invertible, the conditions can be equivalently expressed as the fixed point problem
\begin{equation}
	\label{eq:global minimiser}
	\widehat \rho = \max\{ 0, (U')^{-1} (C - V(x) - W*\widehat\rho) \} \quad \text{in } \mathbb R^d.
\end{equation}
We refer to \cite{BalagueCarrilloLaurentRaoul_Dimensionality,CarrilloDelgadinoPatacchini2019} for a proof of these conditions, as well as for the Euler-Lagrange conditions with respect to other distances between probability measures. Global minimizers in the purely aggregation equation, $U=0$, have also been studied a lot in the last years, see \cite{BLL12,CDM,CCP15,CH17,CFP17,CMMRSV20,F22,DLM22,DLM23,CS23,CS24} for instance. This is a subtopic of the theory of aggregation-diffusion equations on its own with many interesting open problems such as the dimension of the support of global minimizers in certain ranges \cite{BalagueCarrilloLaurentRaoul,BalagueCarrilloLaurentRaoul_Dimensionality}, see the survey \cite{CV15} and the references therein.

\subsection{Steady States}

Steady states are solutions (in a suitable sense) to the time-independent problem
\begin{equation*}
	-\nabla \cdot (\rho \nabla (U'(\rho) + V + W*\rho)) = 0.
\end{equation*}
Under mild assumptions, steady states are also characterised by the densities where the free-energy dissipation \eqref{dissipation} vanishes; that is,
\begin{equation*}
	\int_{\R^d} \rho |\nabla (U'(\rho) + V + W*\rho)|^2 \dx = 0.
\end{equation*}
Partitioning the support of the density, $\{x : \rho > 0\}$, as a union of disjoint open sets $A_i$, the condition is equivalent to finding a density $\rho$ and constants $C_i$ such that
\begin{equation*}
	U'(\rho) + V + W*\rho = C_i \quad \text{in } A_i.
\end{equation*}
These steady states need not be minimisers of the energy.

In particular cases (degenerate diffusion, $W=0$, and certain choices of non-convex $V$, such as a double-well potential), steady states can be created by ``glueing'' disjoint, compactly supported solutions, leading to
\begin{equation*}
	\widehat \rho = (U')^{-1} (C_i - V(x)), \quad \text{in } A_i.
\end{equation*}
These steady states sometimes have a non-trivial basin of attraction.

\subsection{Radial Symmetry vs. Symmetry Breaking}

When $V$ and $W$ are both radially symmetric, we may expect the stationary states of \eqref{aggeqn0} to be radially symmetric as well, yet this is not true in general. For instance, if $W$ is attractive, smooth, and compactly supported, $V$ is zero, and the diffusion is degenerate, there is a single ``bump'' stationary state which is radially symmetric, but multiple copies of the state can be present in a non-interacting way (sufficiently far from each other) to form non-trivial states. In one dimension, the ``gluing'' example of the previous section is valid if $V$ is a symmetric double-well potential. In general, \emph{symmetry breaking} is a well-documented phenomenon in PDEs (e.g. \cite{dolbeault2017SymmetrySymmetryBreaking}), and there are many examples that are derived from particle systems (see for instance \cite{BalagueCarrilloLaurentRaoul_Dimensionality,bertozzi2015RingPatternsTheir}). We highlight \cite{delgadino2022UniquenessNonuniquenessSteady}, where examples of $W$ with infinitely many radially decreasing steady states are constructed.

Even in the case of globally supported interaction potentials, the radial symmetry of steady states is generally difficult to prove. The global minimisers of the free energy can be proven to be radially symmetric when $V$ and $W$ are radially increasing, through the method of Schwarz rearrangement \cite{Talenti2016}. This method replaces $\rho$ by a radially symmetric non-increasing
function $\rho^\star$ for minimizing sequences of the free energy in \eqref{free} such that the measure of the level sets is preserved:
\begin{equation*}
	|\{ x : \rho^\star \ge h \}| = |\{ x: \rho \ge h\}|, \quad \forall h \in \mathbb R.
\end{equation*}
Proving radial symmetry for the steady states of the equation is much more difficult than for the global minimisers of the energy. The idea there is to use continuous Steiner rearrangement (see \cite{Brock2000}), where one continuously ``slides'' each level set. This was used in \cite{Carrillo2019a} to prove radial symmetry of the steady states for $V$ and $W$ radially increasing, and subsequently generalised for potentials which are more singular at the origin \cite{CHMV}.

Establishing whether the radially symmetric and non-increasing steady states are unique (up to translations, if $V=0$) would be useful in order to understand the long-time asymptotics of these models; note, further, that uniqueness of the steady states would also imply uniqueness of the global minimisers for the energy. However, the question of sharp conditions for uniqueness is still open, though important advances have been made under various assumptions on the interaction potentials, see \cite{Lieb83,CGHMV20,CCH21,delgadino2022UniquenessNonuniquenessSteady}. The uniqueness of absolutely continuous stationary states is very important in certain applications, such as neural networks models based on aggregation-diffusion equations, where the uniqueness would provide guarantees on the performance of the network \cite{RV22,BC23}.

\subsection{Steady States Without Aggregation}
When $W = 0$, the condition \eqref{eq:global minimiser} gives an explicit candidate for the global minimiser, provided a value of $C$ consistent with the mass of the steady state can be found (this may not always be possible). Uniqueness follows from a monotonicity argument.

In the case of fast diffusion, steady states might involve Dirac deltas, and their characterisation may require the lower-semi-continuous extension of the free energy,
\begin{equation*}
	\widetilde {\mathcal F}(\mu) = \inf_{\rho_n \rightharpoonup \mu } \liminf_{n} \mathcal F(\rho_n),
\end{equation*}
where $\rho_n$ are densities. Under a radially symmetric potential $V$, the global minimiser in $\mathcal P(\mathbb R^d)$ of $\widetilde {\mathcal F}$ sometimes takes the form
\begin{equation}
	\label{eq:steady state with concentrated part}
	\widehat \mu = (1- \mathfrak m) \delta_0 + \widehat \rho, \quad \text{where }\mathfrak m = \int_{\mathbb R^d} \widehat \rho \dx,
\end{equation}
where $\widehat \rho$ is a density. This is indeed the global attractor in some cases \cite{CarrilloGomezCastroVazquez2022JMPA}, where the density spontaneously forms a Dirac delta at the origin, an effect known as \emph{concentration}.

\subsection{A Clearer Picture for the Homogeneous Setting}
\label{sec:homogeneous case}
A particularly interesting setting is the so-called \emph{homogeneous setting}:
\begin{equation}\label{homog}
	U(\rho) = \frac{\rho^m}{m-1}, \quad W(x) = \chi \frac{|x|^k}{k},
\end{equation}
with $V = 0$ for simplicity, and where we formally denote $|x|^0 / 0 \coloneqq \log |x|$. The parameter $\chi>0$ represents the relative strength of attraction over diffusion in the system. A direct scaling analysis (for $k \ne 0$) leads to
\begin{equation*}
	\mathcal F( \rho_\lambda ) = \lambda^{d(m-1)} \frac{1}{m-1} \int_{\R^d} \rho_1^m(x) \dx + \lambda^{-k} \frac{\chi}{2k} \int_{\R^d} \int_{\R^d} |x-y|^{k} \rho_1(x) \rho_1(y) \dx \dy,
\end{equation*}
where $\rho_\lambda(x):=\lambda^{d} \rho_1(\lambda x )$, for all $\lambda>0$, and where $\rho_1$ is a fixed density. This reveals the critical exponent $m_c = \frac{d-k}{d}$ (the \emph{fair competition range}), where the free energy is homogeneous. The range $m > m_c$ is known as the \emph{diffusion-dominated regime} (where spreading is more energetically favourable than concentration) and $m \in (0,m_c)$ is the \emph{aggregation-dominated regime}, (where the opposite is true). The homogeneous setting has been the subject of significant work, and is relatively well-understood. We summarise the results:
\begin{theorem}[Homogeneous Interaction Potential \cite{CarrilloDelgadinoPatacchini2019,carrillo2019ReverseHardyLittlewood,CCH17}]\label{thmhom}

	\begin{enumerate}
		\item \emph{(Degenerate diffusion - Diffusion dominated)} If $m > 1$ and $k \in ((1-m)d , 0)$, then the energy
		      \begin{equation*}
			      \widetilde{\mathcal F}(\rho) = \begin{dcases}
				      \mathcal F (\rho) & \text{if } \rho \in L^m (\mathbb R^d)                                   \\
				      +\infty           & \text{if } \rho \in \mathcal P(\mathbb R^d) \setminus L^m (\mathbb R^d)
			      \end{dcases}
		      \end{equation*}
		      has a global minimiser for all values of the parameter $\chi>0$.

		\item \emph{(Degenerate diffusion - Fair competition)} If $m > 1$ and $k=(1-m)d$, then there is a dichotomy: the energy $\widetilde{\mathcal F}(\rho)$ has global minimisers if and only if $\chi$ coincides with a critical parameter $\chi_c$. If $0<\chi<\chi_c$, the free energy functional has zero infimum, but it is not achieved. If $\chi>\chi_c$, the free energy functional is not bounded below.

		\item \emph{(Fast diffusion - Aggregation dominated)} If $m \leq \frac{d}{d+k}$ and $k>0$, then the free energy is not bounded below in $\mathcal P(\mathbb R^d) \cap L^\infty (\mathbb R^d)$.

		\item \emph{(Fast diffusion - Diffusion dominated)} If $\frac{d}{d+k}<m<1$ and $k>0$, then the free energy is bounded below in $\mathcal P(\mathbb R^d) \cap L^\infty (\mathbb R^d)$ and there are radially decreasing global minimisers of a suitably extended free energy functional (allowing for a possible partial mass concentration at the origin).
	\end{enumerate}
\end{theorem}

The proofs of these results hinge drastically on the sharp constants of related functional inequalities.

\subsection{Degenerate Diffusion}
The case of degenerate diffusions is handled with the classical Hardy-Littlewood-Sobolev (HLS) inequality; to be precise, with the following variant obtained by interpolation. Let $k \in (-d,0)$, and $m \in (1,2)$ such that $d(m-1)+k=0$. Then there exists a positive and finite optimal constant $C_*$ such that, for $f \in L^1(\R^d) \cap L^m(\R^d)$, we have
\begin{equation}\label{vhls}
	\left|\iint_{\R^d\times \R^d} f(x){|x-y|^{k}}f(y) \dx\dy\right| \leq C_* ||f||^{\frac{d+k}{d}}_1 ||f||^m_m.
\end{equation}
This inequality permits the comparison of competing terms in the free energy $\mathcal{F}[\rho]$ in \eqref{free}, where $V=0$ and $U$ and $W$ are given by \eqref{homog}, since the internal energy is positive for $m>1$, and the interaction energy is negative.
By defining the critical interaction strength $\chi_c$ as the biggest $\chi$ such that $\mathcal{F}[\rho] \geq 0$ for all $\rho$, one can establish the following:

\begin{enumerate}
	\item \emph{(Critical case $\chi=\chi_c$)} The variant of the HLS inequality is equivalent to $\mathcal F[\rho] \geq 0$. Stationary states are equivalent to optimisers of the variant of the HLS inequality. Moreover, for the critical case, all of them are global minimisers of the free energy. Furthermore, there exists $\rho^* \in \mathcal P(\mathbb R^d)\cap L^m(\R^d)$, non-increasing, radially symmetric, and with $||\rho_*||_m=1$ such that $\mathcal F[\rho_*]=0$; all steady states are translations and dilations of $\rho_*$. These are compactly supported, bounded, and integrable enough for the singular convolution to make sense in the full range. Moreover, at the boundary, they have limited regularity, like the Barenblatt profiles of the porous medium equation \cite{Va}.

	\item \emph{(Subcritical case $0<\chi<\chi_c$, $k<0$)} There are stationary states in rescaled coordinates, which lead to self-similar solutions in original coordinates.

	\item \emph{(Supercritical case $\chi>\chi_c$, $k<0$)} The free energy $\mathcal F[\rho]$ is not bounded below. Moreover, global solutions that decay fast at infinity do not exist; this is shown by a contradiction argument, as for the classical Keller-Segel model \cite{DP,BlanchetDolbeaultPerthame}.

\end{enumerate}

In short, the variant of the HLS inequality \eqref{vhls} for the degenerate diffusion range $m>1$ and the fair competition regimes $m=1-\frac{k}d$ leads to the same dichotomy as for the classical Keller-Segel model \cite{DP, bdp} (which uses the logarithmic Hardy-Littlewood-Sobolev inequality \cite{CL}). For more details, we refer for more details to \cite{DP,BlanchetDolbeaultPerthame,BCL,BlanchetCarrilloMasmoudi,BCC,CCH17,CHMV}.

\subsection{Fast Diffusion}

\begin{figure}\label{fig:fastdif}
	\centering
	\begin{tikzpicture}
		\draw[->] (-0.2,0) -- (10,0) node[right] {$k$};
		\draw[->] (0,-0.2) -- (0,3.5) node[above] {$m$};

		\draw[scale=2.5,fill=gray!50!white] plot[smooth,samples=100,domain=0:0.43] ({\x},{1/(1+\x)}) --
		plot[smooth,samples=100,domain=0.43:0] (\x,{1});
		\draw[draw=gray!50!white,fill=gray!50!white]
		plot[smooth,samples=100,domain=1:10] (\x,{1.75}) --
		plot[smooth,samples=100,domain=10:1] (\x,{2.5});

		\draw[scale=2.5,fill=gray!50!white] plot[smooth,samples=100,domain=0:4] ({\x},{3.5/(3.5+\x)}) --
		plot[smooth,samples=100,domain=4:0] (\x,{1.});

		\draw[scale=2.5,domain=0:1.2,smooth,variable=\x,red,thick] plot ({\x},{1-\x}) node[right] {$m=1-\frac{k}{d}$};
		\draw[scale=2.5,domain=0:4,smooth,variable=\x,black,thick] plot ({\x},{1/(1+\x)}) node[right] {$m=\frac{d}{d+k}$};
		\draw[scale=2.5,domain=0:4,smooth,variable=\x,red, dashed] plot ({\x},{3.5/(3.5+\x)}) node[right] {$m=\frac{2d}{2d+k}$};
		\draw[thick](0,2.5) -- (10,2.5) node[right] {$m=1$};
		\draw[thick,blue] (0,1.75) -- (10,1.75) node[right] {$m=\frac{d-2}d$};

		\draw (2.4,-0.2) node {$d$};
		\draw (6.5,0.35) node {{\large \bf I}};
		\draw (6.5,1.5) node {{\large \bf II}};
		\draw (6.5,3.3) node {{\large \bf III}};

		\draw[->,dashed] (2,0) node[below] {$4$} -- (2,1.75);
		\draw[->,thick,red] (2,1.4) -- (2,1.75);
	\end{tikzpicture}
	\caption{Sketch of the main results known for the fast diffusion range. The dark grey region corresponds to the area in which global minimisers of the free energy are integrable. The arrowed red line is the region for which we know the exact value of $m$ separating integrable from partially mass-concentrated global minimisers. The remaining white region in Zone II corresponds to the set of parameters in which the separation of integrable and partial mass concentrated global minimisers is not known yet.
	}
\end{figure}
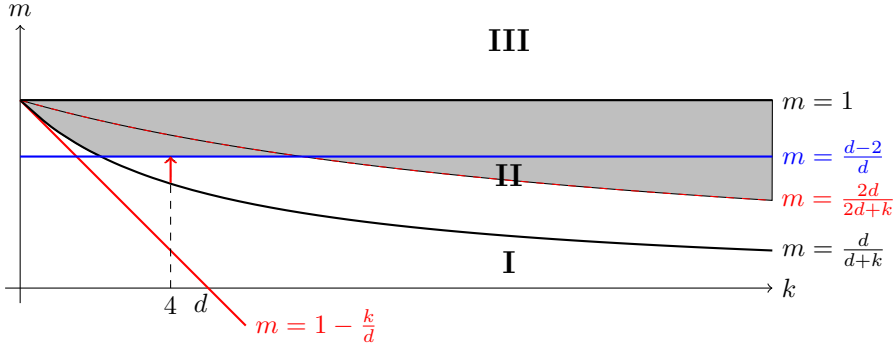

In \cite{carrillo2019ReverseHardyLittlewood}, the authors introduce a reversed HLS (RHLS) inequality to deal with the diffusion-dominated regime:
\begin{equation}
	\label{eq:reversed HLS}
	\iint_{\mathbb R^d \times \mathbb R^d} f(x) |x-y|^k f(y) \dx \dy \geq C \left( \int_{\mathbb R^d} f(x) \dx \right)^{\alpha} \left( \int_{\mathbb R^d} f(x)^m \dx \right)^{\frac{2-\alpha}{m}} .
\end{equation}
Again, the RHLS inequality serves to control the negative part of the free energy $\mathcal{F}[\rho]$ in \eqref{free}, where $V=0$ and $U$ and $W$ are given by \eqref{homog}. To be precise, we can control the negative part (the internal energy in the fast diffusion range), by the positive part (the interaction energy, which now is a sort of non-local moment). Notice that the roles of the competing mechanisms in the free energy functional has swapped with respect to the degenerate diffusion case. The interesting fact for $m<1$ is that this can only happen for $m>\frac{d}{d+k}$, as stated in \cref{thmhom}. Moreover, the homogeneity line $m=1-\frac{k}d$ does not play any role, since it lies below the curve $m=\frac{d}{d+k}$ (and is tangent to it at $m=1$); thus, it lies in the region where the free energy is not bounded below. We refer to \cref{fig:fastdif} for a sketch of the main results.

We denote by $L^1\cap L^m(\R^d)$ the set of integrable functions on $\R^d$ such that $\rho^m$ is an integrable function (recall $0<m<1$ in this section, so this is not a Banach space). The main results can be summarised as:

\

\begin{itemize}
	\item Zone I, $m\leq \frac{d}{d+k}$: The free energy in $L^1\cap L^m(\R^d)$ is not bounded below up to $m=\frac{d}{d+k}$. The dilation scaling relation, $m=1-\frac{k}d$, does not play a role.
	      \\
	\item Zone II, $\frac{d}{d+k}<m<1$: The reversed HLS inequality \eqref{eq:reversed HLS} holds.
	      \begin{enumerate}
		      \item The minimiser is attained in $L^1\cap L^m(\R^d)$ in the dark grey region. It is a positive and radially decreasing density (up to translations).

		      \item In the rest of this zone, the minimiser is attained by suitably extending the free energy functional to the set of probability measures. It consists of a radially decreasing function in $\mathrm L^1\cap\mathrm L^m(\R^d)$ with {\bf possibly} a Dirac Delta at the origin (modulo translations). We call this a possible partial mass concentration steady state \eqref{eq:steady state with concentrated part}.

		      \item The black-red dashed line corresponds to the only case of the RHLS inequality known previous to \cite{carrillo2019ReverseHardyLittlewood}, see \cite{DZ16,NN17}.

		      \item The arrow-red line is the area where partial mass concentration happens if $m$ is below an explicit sharp value of the exponent $m$ as proven in \cite{CDFL}.
	      \end{enumerate}
	      \phantom{dddd}
	\item Zone III, $m\geq 1$: The minimiser is attained and it is a radially symmetric, decreasing, and compactly supported density, modulo translations.
\end{itemize}

\

Part of these results can be found in \cite{carrillo2019ReverseHardyLittlewood}. They also appeared as a preprint \cite{CarrilloDelgadino2018} written in terms of \eqref{aggeqn0} and its free energy functional \eqref{free}, instead of functional inequalities. The authors also give suitable Euler-Lagrange equations for $\widehat \mu$, which are more involved than those in \Cref{sec:Euler Lagrange} due to the presence of the Dirac delta. It is worth remarking that $k = 2$ is trivial, since $W*\rho$ is explicit for radial functions. The case $k = 4$ was studied in full detail in \cite{CDFL}. There, the authors characterise precisely the dimensions $d$ and indices $m$ for the global minimiser to have a Dirac delta with partial mass at the origin, i.e., \eqref{eq:steady state with concentrated part} with $\mathfrak m = \|\widehat \rho\|_{L^1} < 1$.

There are a number of interesting open problems in this direction. Finding the general relation between $m$ and $k$ for which there is concentration in the white region of zone II is a very interesting open question. Studying the evolution problem associated to the gradient flow of these free energy functionals is also a difficult problem. Characterising the conditions on the initial data leading to a partial mass concentration at the origin, either asymptotically, or in finite time, eventually converging to the final global minimiser of the energy with partial mass concentrated at the origin, is an even more challenging open problem.

\subsection{Bifurcation of Steady States}

The existence of multiple solutions of the fixed-point problem \eqref{eq:global minimiser} is particularly interesting as we increase the influence of $W$, for instance, by including a multiplicative parameter $\chi>0$, as in $W(x) = \chi \frac{|x|^k}{k}$. A simple example is provided by the case of linear diffusion $U = \rho \log \rho$ and $V = 0$; there, we can re-state the problem as finding all solutions of the non-linear integral equation:
\begin{equation*}
	\widehat \rho (x) = \frac{e^{-W*\widehat \rho (x) }}{\displaystyle\int_{\mathbb R^d} e^{-W*\widehat \rho (y)} \dy}, \quad \text{for all } x \in \mathbb R^d.
\end{equation*}
\cite{BCCD16,carrillo2020LongTimeBehaviourPhase} are two examples of works dealing with the bifurcation branches for \eqref{aggeqn0}: the first one, on the whole space (reducible to a one dimensional bifurcation problem), and the second one, with periodic conditions (by using Crandall-Rabinowitz theory). The main results in this direction show the the characterisation of bifurcation branches and the multiplicity of solutions for this family of equations is a complex and intricate process. The bifurcations from constant solutions in the periodic case or from symmetric solutions in the whole space have been reported for both linear and non-linear diffusions \cite{BCCD16,BCDPZ18,CCWWZ20,CG21}.

In the periodic setting an almost full characterization of all bifurcation points from constants steady states can be given in terms of Fourier modes of the interaction potential $W$, see \cite{carrillo2020LongTimeBehaviourPhase,CG21}. The classical case is the Kuramoto model for synchronization, which leads to a single branch of bumped solutions at a well-know critical value \cite{kuramoto1975international,SM91,acebron2005kuramoto,CCP19}.
The stability analysis of these stationary solutions for \eqref{aggeqn0} is only known in very special situations, and a full clarification of the basin of attractions or stable/unstable manifolds of bifurcation branches is still out of reach. Some of these aspects need a full exploration with suitable numerical methods, see \cite{BCH2020}.

To conclude, we discuss another interesting aspect of this problem: the continuity of the behaviour of the system when a bifurcation takes place. Could it happen, in the periodic case, that there are branches of global minimisers different from the constant solution for values of the parameter $\chi$ before the first bifurcation branches at parameter $\chi_*$? This is known as a \emph{discontinuous phase transition}. On the other hand, the classical Kuramoto model is the archetypical example of continuous phase transition, in which a branch of stable equilibria appears at a critical value in a continuous way emanating from the constant steady state. Sufficient conditions for continuous and discontinuous phase transitions to occur, both in the linear and non-linear diffusion cases, have been given in \cite{carrillo2020LongTimeBehaviourPhase,CG21}.

\section{Numerical Schemes}

\begin{figure}
	\centering

	\subfloat[$t=0$]{\includegraphics[width=0.32\textwidth,trim={0.2cm 0.4cm 0.2cm 0.85cm},clip]{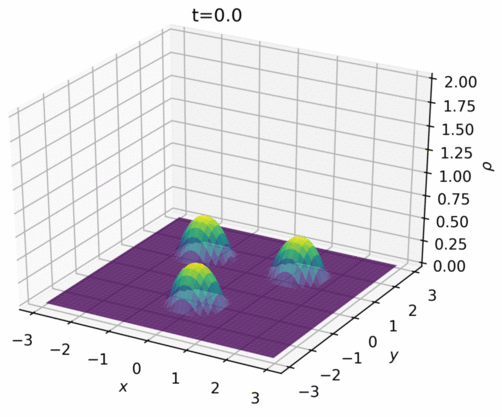}}
	\subfloat[$t=155$]{\includegraphics[width=0.32\textwidth,trim={0.2cm 0.4cm 0.2cm 0.85cm},clip]{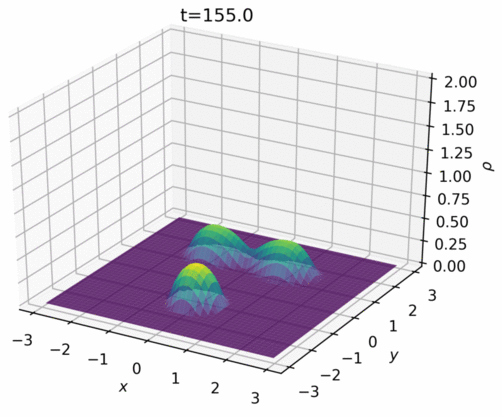}}
	\subfloat[$t=200$]{\includegraphics[width=0.32\textwidth,trim={0.2cm 0.4cm 0.2cm 0.85cm},clip]{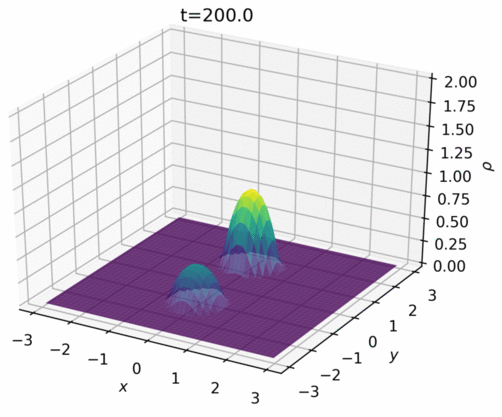}}

	\subfloat[$t=275$]{\includegraphics[width=0.32\textwidth,trim={0.2cm 0.4cm 0.2cm 0.85cm},clip]{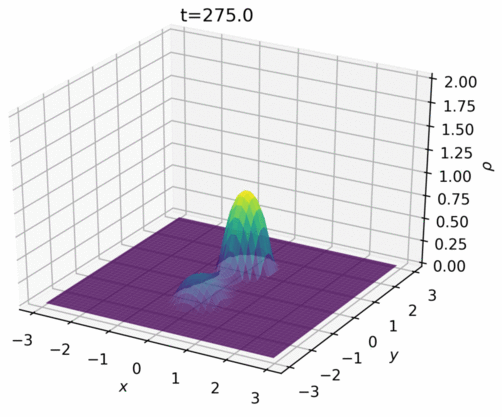}}
	\subfloat[$t=360$]{\includegraphics[width=0.32\textwidth,trim={0.2cm 0.4cm 0.2cm 0.85cm},clip]{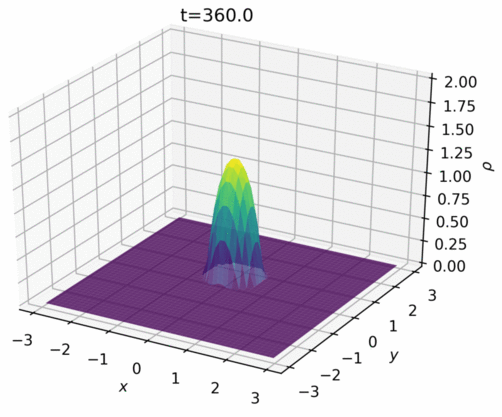}}
	\subfloat[Relative free energy]{\includegraphics[width=0.32\textwidth,trim={3cm 1cm 6cm 6cm},clip]{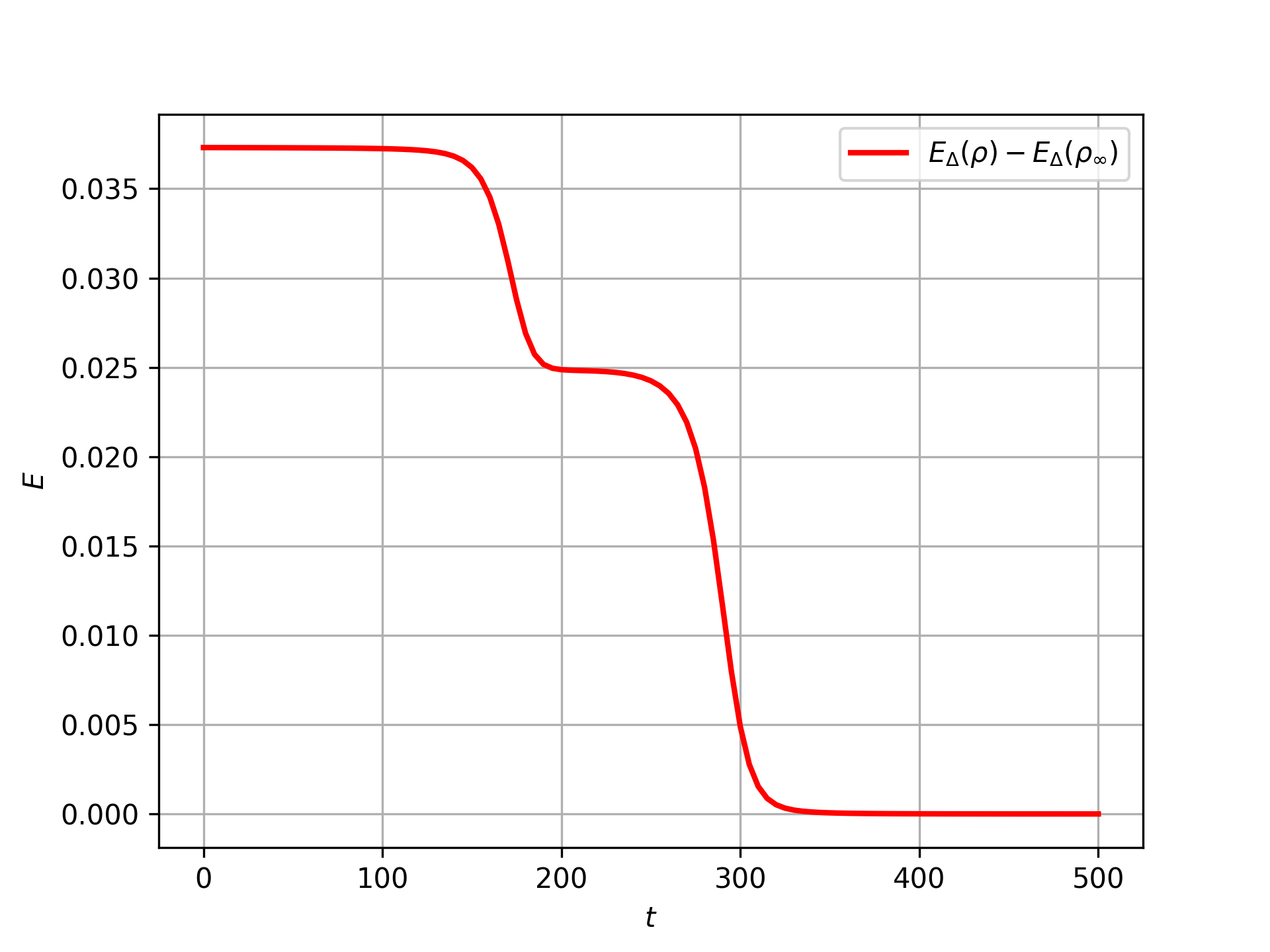}}

	\caption{Metastability in the aggregation-diffusion equation \eqref{aggeqn0} under a force with exponential decay from \cite{BCH2020}. Long periods of slow motion are interspersed by rapid merging phases until a steady state is reached.}
	\label{fig:bumps}
\end{figure}

\begin{figure}
	\centering

	\subfloat[Cell adhesion]{\includegraphics[width=0.49\textwidth,trim={0cm 3cm 0cm 0cm},clip]{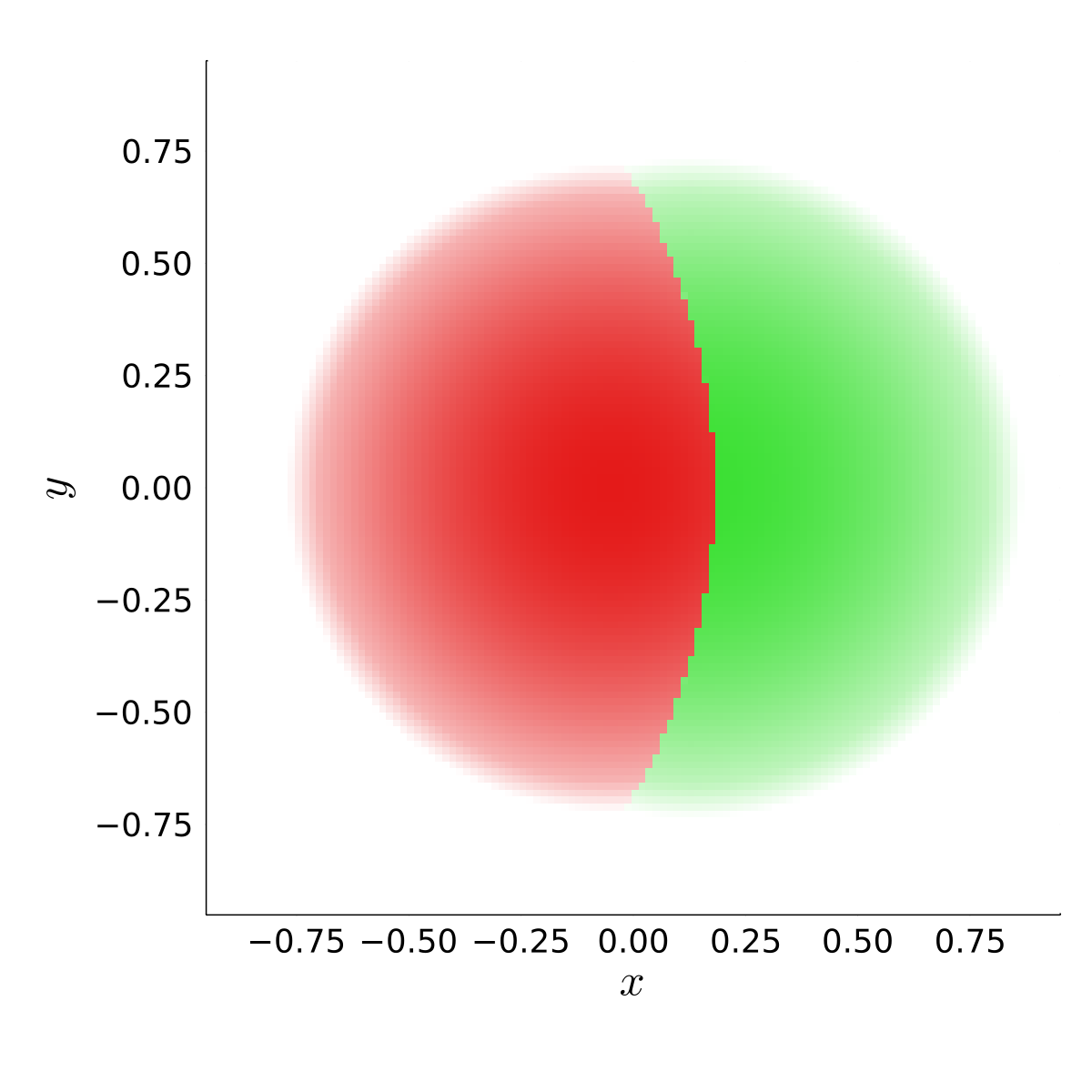}\label{fig:mastercard_a}}
	\subfloat[Adhesion with diffuse boundary]{\includegraphics[width=0.49\textwidth,trim={0cm 3cm 0cm 0cm},clip]{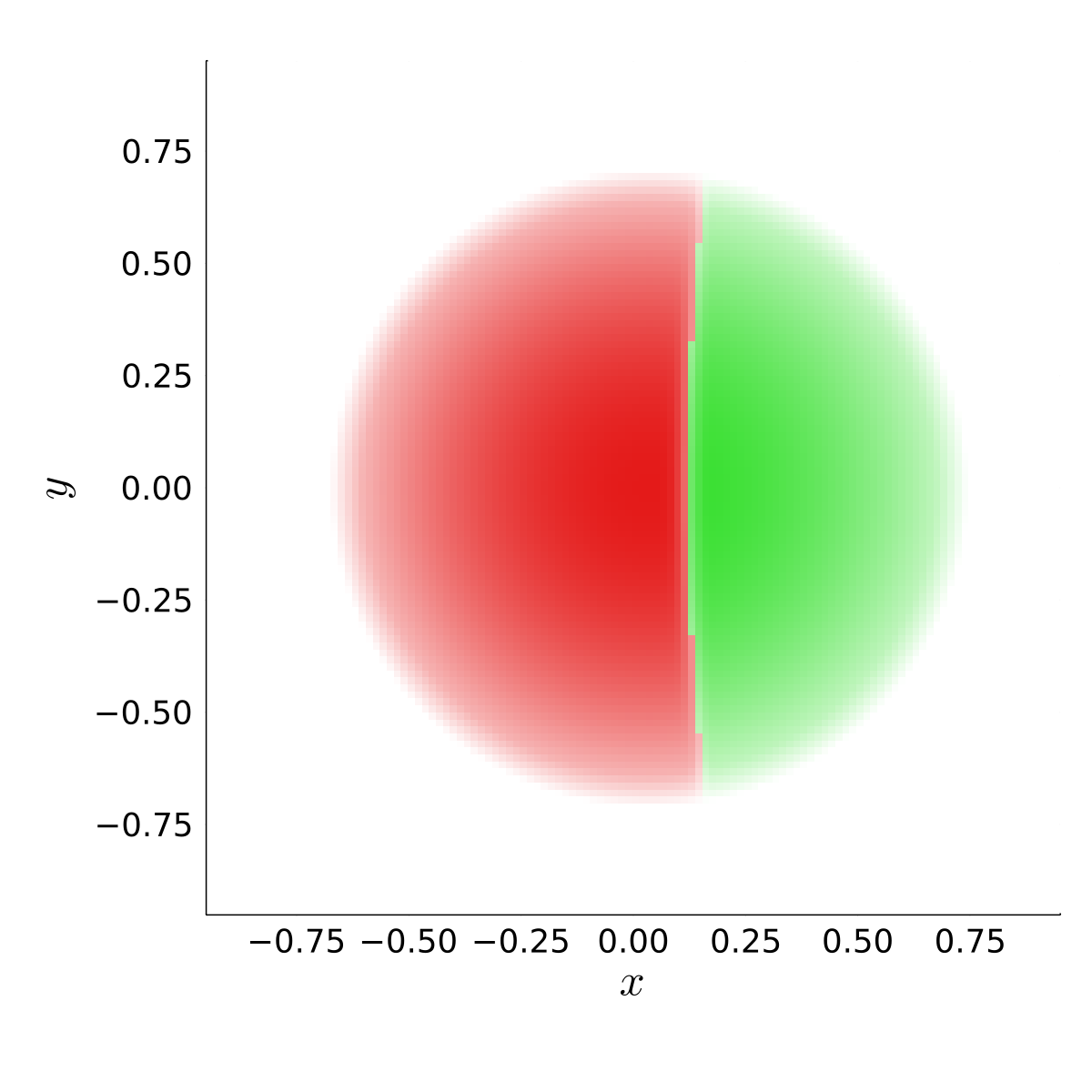}\label{fig:mastercard_b}}

	\subfloat[Cellular engulfment]{\includegraphics[width=0.49\textwidth,trim={0cm 3cm 0cm 0cm},clip]{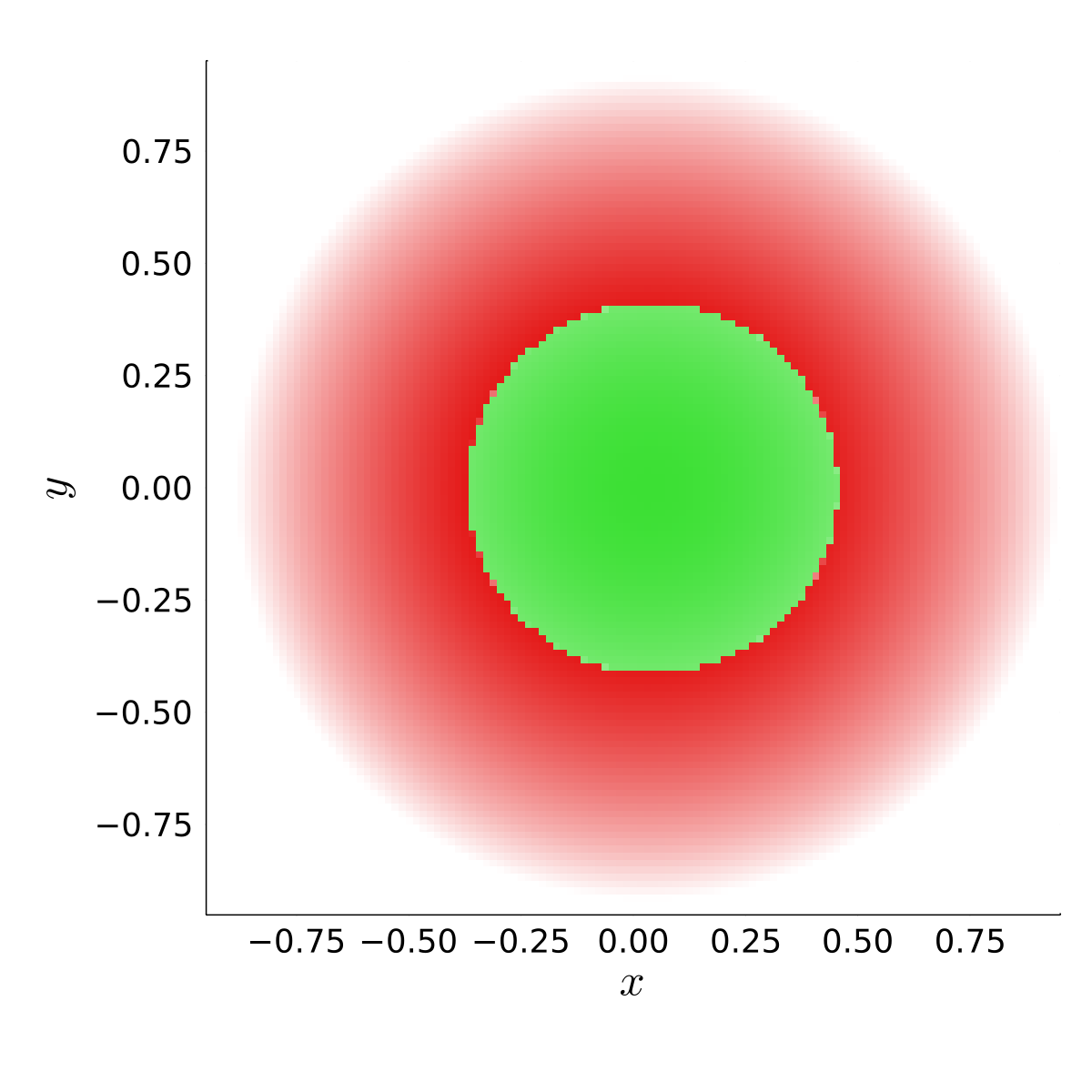}\label{fig:mastercard_c}}
	\subfloat[Engulfment limited by volume constraints]{\includegraphics[width=0.49\textwidth,trim={0cm 3cm 0cm 0cm},clip]{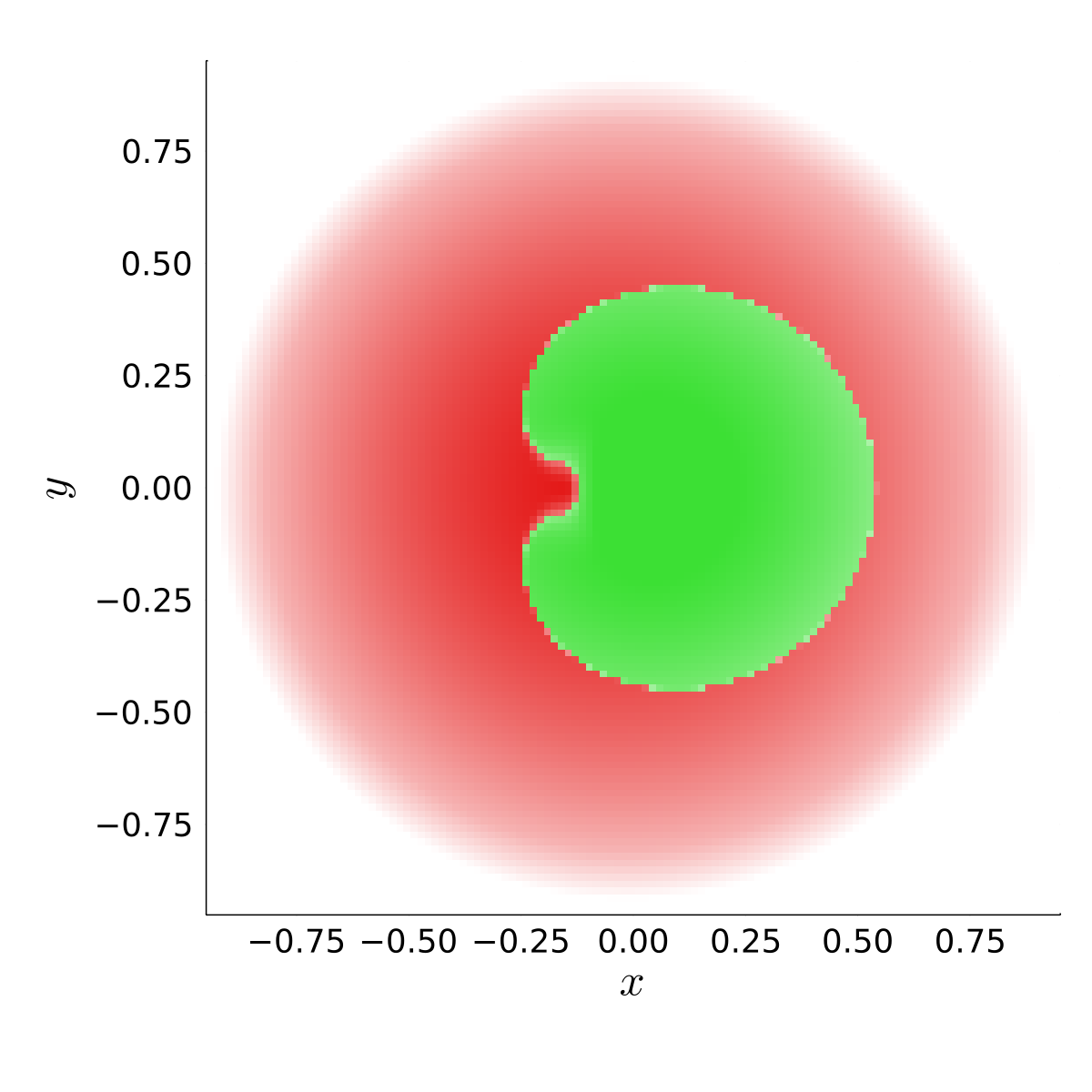}\label{fig:mastercard_d}}

	\caption{Cell aggregates mediated by cellular adhesion in model \eqref{model_chs} under different parameters, coloured by population, adapted from \cite{BCH2023}.}
	\label{fig:mastercard}
\end{figure}

The aggregation-diffusion equation \eqref{aggeqn0} can be solved numerically in myriad ways, such as mesh-based schemes \cite{CHM2023,JPZ2024}, particle methods \cite{CB2015,CCC2018,CCP2019}, and methods based in the JKO iteration \cite{CCW2021}. In some cases, classical numerical methods for the Fokker-Planck equation such as the Scharfetter-Gummel scheme have been generalised to \eqref{aggeqn0} \cite{SS2021,HST2023}. Other approaches treat the equation as a conservation law in order to establish structure-preserving properties \cite{SCS2018,ABP2019,CCF2020} (as was done in \cite{BF2012} for non-linear diffusions).

An example of the latter is \cite{CCH2015}, which introduces the finite-volume spatial semi-discretisation (presented in one dimension for simplicity)
\begin{subequations}\label{eq:scheme_continuous}
	\begin{align}
		 & \der{\rho\i}{t} + \frac{F\ih - F\imh}{\Dx} = 0,
		 &                                                 & F\ih = \rho\i \pos{u\ih} + \rho\ip \neg{u\ih}, \\
		 & u\ih=-\frac{\xi\ip-\xi\i}{\Dx},
		 &                                                 & \xi\i=H'(\rho\i)+V\i+(W\ast\rho)\i,
	\end{align}
\end{subequations}
where $\pos{s}\coloneqq \max\set{0,s}$ and $\neg{s}\coloneqq \min\set{0,s}$, and where the index $i$ represents the value at a point $x_i$ in a discretised domain. The terms $V\i$ and $(W\ast\rho)\i$ are suitable approximations of the potential $V(x)$ and the convolution term $\prt*{W\ast\rho}(x)$. A fully discrete scheme was given in \cite{BCH2020,BCM2020}:
\begin{subequations}\label{eq:scheme}
	\begin{align}
		 & \frac{\rho\i\np - \rho\i\n}{\Dt} + \frac{F\ih\np - F\imh\np}{\Dx} = 0,
		 &                                                                        & F\ih\np = \rho\i\np \pos{u\ih\np} + \rho\ip\np \neg{u\ih\np}, \\
		 & u\ih\np=-\frac{\xi\ip\np-\xi\i\np}{\Dx},
		 &                                                                        & \xi\i\np=H'(\rho\i\np)+V\i+\prt{W\ast\rho\nss}\i,
	\end{align}
\end{subequations}
where $\rho\i\nss = (\rho\i\n + \rho\i\np)/2$, and where the index $n$ indicates the value at time $t\n = n\Dt$. Schemes \eqref{eq:scheme_continuous} and \eqref{eq:scheme} capture the (respectively) semi-discrete and fully discrete analogues of three crucial properties of the aggregation-diffusion equation \eqref{aggeqn0}: the conservation of mass, the non-negativity of solutions, and the dissipation of the free energy.

The conservation of mass and non-negativity properties of the continuous problem can be rigorously established by studying \eqref{aggeqn0} as a conservation law, but they are intuitively suggested by the derivation of the equation as a mean-field limit: the equation models a density of individuals (which cannot be negative), and their motion in space, rather than their birth or death (thus the number of individuals is constant in time). Because these principles are so intrinsic to the model, it is desirable to preserve them when conducting simulations, which is the case in these numerical schemes:

\begin{proposition}[Conservation of mass and non-negativity]
	Let $\rho\i(t)$ be a solution to \eqref{eq:scheme_continuous} with $\rho\i(0)\geq 0$ and $\norm{\rho\i(0)}_{L^1} = M$. Then, $\rho\i(t)\geq 0$ and $\norm{\rho\i(t)}_{L^1} = M$ for all $t>0$.

	Similarly, let $\rho\i\n$ be a solution to \eqref{eq:scheme} with $\rho\i^0\geq 0$ and $\norm{\rho\i^0}_{L^1} = M$. Then, $\rho\i\n\geq 0$ and $\norm{\rho\i\n}_{L^1} = M$ for all $n>0$.
\end{proposition}

The dissipation of the free energy in the aggregation-diffusion equation was already discussed in \cref{sec:gradient_flows}. While solutions to \eqref{aggeqn0} can be rigorously constructed as gradient flows of the energy \eqref{free}, the milder energy dissipation identity \eqref{dissipation} can be shown formally by directly computing the time derivative of the free energy along solutions to the equation assuming suitable boundary conditions. The last term in the identity reveals that the non-negativity of the solution $\rho$ is in fact a crucial ingredient of the energy dissipation. In part because said non-negativity is preserved by the numerical schemes above, discrete energy dissipation \emph{inequalities} can be established:
\begin{theorem}[Energy dissipation]\label{th:dissipation}
	Let $\rho\i(t)$ be a solution to \eqref{eq:scheme_continuous} with $\rho\i(0)\geq 0$. Define the \emph{semi-discrete free energy} by
	\begin{equation}
		\mathcal{F}_\Delta(t) =
		\sum\i H\prt{\rho\i(t)} \Dx
		+\sum\i V\i\rho\i(t) \Dx
		+\frac{1}{2}\sum\i\sum\k \Wik \rho\i(t) \rho\k(t) \Dx^2.
	\end{equation}
	Then,
	\begin{equation}\label{eq:dissipation_semi}
		\der{}{t} \mathcal{F}_\Delta(t)
		\leq
		-\sum\i \min\set*{\rho\i(t), \rho\ip(t)}\abs{u\ih}^2 \Dx \leq 0.
	\end{equation}

	Similarly, let $\rho\i\n$ be a solution to \eqref{eq:scheme} with $\rho\i^0\geq 0$. Define the \emph{fully discrete free energy} by
	\begin{equation}
		\mathcal{F}_\Delta\n =
		\sum\i H\prt{\rho\i\n} \Dx
		+\sum\i V\i\rho\i\n \Dx
		+\frac{1}{2}\sum\i\sum\k \Wik \rho\i\n \rho\k\n \Dx^2.
	\end{equation}
	Then,
	\begin{equation}\label{eq:dissipation_fully}
		\frac{\mathcal{F}_\Delta\np - \mathcal{F}_\Delta\n}{\Dt}
		\leq
		-\sum\i \min\set*{\rho\i\np, \rho\ip\np}\abs{u\ih\np}^2 \Dx \leq 0.
	\end{equation}
\end{theorem}

Recalling the definition of the velocity field $u=-\nabla\xi$ in the aggregation-diffusion equation \eqref{aggeqn0}, the inequalities \eqref{eq:dissipation_semi} and \eqref{eq:dissipation_fully} can be understood as discretisations of the dissipation identity \eqref{dissipation}.

These properties do, in fact, generalise quite well. Scheme \eqref{eq:scheme_continuous} is also given with a second-order spatial discretisation (using flux limiters) and in two dimensions in \cite{CCH2015}, and can be trivially generalised to higher dimensional settings, always preserving the conservation of mass, non-negativity, and energy dissipation properties. The properties do not survive a trivial discretisation in time (for instance, a fully explicit treatment), and their preservation may require stringent restrictions on the time step $\Dt$. However, the time discretisation \eqref{eq:scheme} preserves the properties unconditionally, and is also generalised to higher dimensions in \cite{BCH2020}. The latter introduces a new type of dimensional splitting for \eqref{eq:scheme} as well which nevertheless preserves the dissipation \cref{th:dissipation}, even in the presence of general non-local potentials $W$; this is needed in practice, as the implicit treatment of higher-dimensional problems can be very costly.

\Cref{fig:bumps} reproduces a numerical example from \cite{BCH2020}: a solution of the aggregation-diffusion equation \eqref{aggeqn0} which exhibits metastability. The particles undergo non-linear diffusion and experience a global but exponentially decaying attraction force. As a result, the density can be organised into compactly supported ``bumps'', each approximately a steady state of \eqref{aggeqn0}, which interact with each other very weakly. These bumps approach each other very slowly over large periods of time, before rapidly merging into a new bump. The slow and fast dynamics are best understood by looking at the free energy and its plateaus. Structure preserving schemes such as \eqref{eq:scheme_continuous} are useful to simulate such a scenario because a large time step can be used without introducing stability issues, oscillations around the boundary of the compact support of each bump, or numerical artifacts in the evolution of the free energy.

The techniques of \cite{CCH2015,BCH2020} have also been generalised to equations with non-linear mobilities, as well as systems of aggregation-diffusion equations in \cite{BCH2023}, always preserving non-negativity and dissipation properties, as well the upper bounds that may be introduced by the mobility terms. Their setting includes the mean-field model \eqref{model_chs} discussed above, and preserves the energy dissipation in the setting where the equation itself is a gradient flow. These ideas have also been extended to other types of gradient flows, such as the Cahn-Hilliard equation \cite{BCK2023}, a gradient flow in $H^{-1}$.

\Cref{fig:mastercard} adapts a numerical test of \cite{BCH2023}, which performs cell adhesion experiments similar to those of \cref{fig:figtethumbnail}. The cellular model \eqref{model_chs} is solved with initially separate populations: just as in \cite{carrillo2019population}, the experimental behaviour described in \cite{katsunuma2016synergistic} can be reproduced by adjusting the adhesion parameters between the cells. Full mixing (not shown in the figure) can be achieved by making all adhesion strengths equal. \Cref{fig:mastercard_a,fig:mastercard_b} show scenarios where the populations remain strongly cohesive, but also attract each other via weaker adhesion; depending on the specific strengths, a clean or a diffusive boundary between the populations can be achieved. \Cref{fig:mastercard_c} shows a scenario where one population completely envelops the other, forming the halo described in \cref{sec:math_bio}. The same scenario with a non-linear mobility function is shown in \cref{fig:mastercard_d}; the additional term models volume constraint effects, where the combined cell density cannot exceed a threshold. This limit forces cells to move around the areas of highest density (such as the point where the two populations initially come into contact, which forms a pinch in the figure), leading to asymmetric behaviour.

\section{Conclusions}

Our goal in this exposition has been to showcase some of the most interesting aspects in the study of aggregation-diffusion equations and systems developed over the past 30 years. The field continues to grow, and while many questions have received answers in recent years, many remain open. To recognise our limited understanding, we remark that long-time asymptotics, bifurcation branches, stable and unstable manifolds, and basins of attraction are only understood for a handful of these models. The applications of aggregation-diffusion equations are vast: mathematical biology, data science, machine learning, and optimisation, to name but a few. The opportunities for applied mathematicians to continue exploring this fascinating topic are boundless.

\bibliographystyle{abbrv}

\bibliography{references,theory,AggDiffBookChapterBiblioFile}

\end{document}